\documentclass[12pt]{article}%
\usepackage{graphicx}
\usepackage[intlimits]{amsmath}
\usepackage{latexsym}
\usepackage{amsfonts}
\usepackage{amssymb}%
\makeatletter
\newfont{\footsc}{cmcsc10 at 8truept}
\newfont{\footbf}{cmbx10 at 8truept}
\newfont{\footrm}{cmr10 at 10truept}
\newtheorem{theorem}{Theorem}

\newtheorem{conjecture}[theorem]{Conjecture}

\newtheorem{problem}[theorem]{Problem}
\newtheorem{proposition}[theorem]{Proposition}

\newenvironment{proof}[1][Proof]{\noindent{\textbf {#1}  }}  {\hfill$\Box$\bigskip}

\begin{document}

\title{Eigenvalue problems of Nordhaus-Gaddum type}
\author{Vladimir Nikiforov\\Department of Mathematical Sciences, University of Memphis, \\Memphis TN 38152, USA}
\maketitle

\begin{abstract}
Let $G$ be a graph with $n$ vertices and $m$ edges and let $\mu_{1}\left(
G\right)  \geq...\geq\mu_{n}\left(  G\right)  $ be the eigenvalues of its
adjacency matrix. We discuss the following general problem. For $k$ fixed and
$n$ large, find or estimate%
\[
f_{k}\left(  n\right)  =\max_{v\left(  G\right)  =n}\text{ }\left\vert \mu
_{k}\left(  G\right)  \right\vert +\left\vert \mu_{k}\left(  \overline
{G}\right)  \right\vert .
\]
In particular we prove that
\[
\frac{4}{3}n-2\leq f_{1}\left(  n\right)  <\left(  \sqrt{2}-c\right)  n
\]
for some $c>8\times10^{-7}$ independent of $n.$ We also show that
\begin{align*}
\frac{\sqrt{2}}{2}n-3  &  <f_{2}\left(  n\right)  <\frac{\sqrt{2}}{2}n,\\
\frac{\sqrt{2}}{2}n-3  &  <f_{n}\left(  n\right)  \leq\frac{\sqrt{3}}{2}n.
\end{align*}


\end{abstract}

\section{Introduction}

Our notation is standard (e.g., see \cite{Bol98}, \cite{CDS80}, and
\cite{HoJo88}); in particular, all graphs are defined on the vertex set
$\left\{  1,2,...,n\right\}  =\left[  n\right]  $ and $G\left(  n,m\right)  $
stands for a graph with $n$ vertices and $m$ edges. We write $\Gamma\left(
u\right)  $ for the set of neighbors of the vertex $u$ and set $d\left(
u\right)  =\left\vert \Gamma\left(  u\right)  \right\vert .$ Given a graph $G$
of order $n,$ we assume that the eigenvalues of the adjacency matrix of $G$
are ordered as $\mu\left(  G\right)  =\mu_{1}\left(  G\right)  \geq...\geq
\mu_{n}\left(  G\right)  $. As usual, $\overline{G}$ denotes the complement of
a graph $G$ and $\omega(G)$ stands for the clique number of $G.$

Nosal \cite{Nos70} showed that for every graph $G$ of order $n,$
\begin{equation}
n-1\leq\mu\left(  G\right)  +\mu\left(  \overline{G}\right)  <\sqrt{2}n.
\label{Nosin}%
\end{equation}
Quite of attention has been given to second of these inequalities. In
\cite{Nik02} it was shown that%
\begin{equation}
\mu\left(  G\right)  +\mu\left(  \overline{G}\right)  \leq\sqrt{\left(
2-\frac{1}{\omega(G)}-\frac{1}{\omega(\overline{G})}\right)  n\left(
n-1\right)  }, \label{Nikin}%
\end{equation}
improving earlier results in \cite{Hon95}, \cite{HoSh00}, \cite{Li96}, and
\cite{Zho97}. Unfortunately inequality (\ref{Nikin}) is not much better then
(\ref{Nosin}) when both $\omega(G)$ and $\omega(\overline{G})$ are large
enough. Thus, it is natural to ask whether $\sqrt{2}$ in (\ref{Nosin}) can be
replaced by a smaller absolute constant for $n$ sufficiently large. In this
note we answer this question in the positive but first we state a more general problem.

\begin{problem}
For every $1\leq k\leq n$ find%
\[
f_{k}\left(  n\right)  =\max_{v\left(  G\right)  =n}\left\vert \mu_{k}\left(
G\right)  \right\vert +\left\vert \mu_{k}\left(  \overline{G}\right)
\right\vert .
\]

\end{problem}

It is difficult to determine precisely $f_{k}\left(  n\right)  $ for every $n$
and $k,$ so at this stage it seems more practical to estimate it
asymptotically. In this note we show that
\begin{equation}
\frac{4}{3}n-2\leq f_{1}\left(  n\right)  <\left(  \sqrt{2}-c\right)  n
\label{mainin1}%
\end{equation}
for some $c>8\times10^{-7}$ independent of $n.$ For $f_{2}\left(  n\right)  $
we give the following tight bounds%
\begin{equation}
\frac{\sqrt{2}}{2}n-3<f_{2}\left(  n\right)  <\frac{\sqrt{2}}{2}n.
\label{mainin2}%
\end{equation}

We also show that
\begin{equation}
\frac{\sqrt{2}}{2}n-3<f_{n}\left(  n\right)  \leq\frac{\sqrt{3}}{2}n
\label{mainin3}%
\end{equation}

Finally for fixed $k,$ $2<k<n,$ and $n$ large, we prove that%
\begin{align*}
\left\lfloor \frac{n}{k}\right\rfloor -1  &  \leq f_{k}\left(  n\right)
\leq\sqrt{\frac{2}{k}}n,\\
\left\lfloor \frac{n}{k}\right\rfloor +1  &  \leq f_{n-k}\left(  n\right)
\leq\sqrt{\frac{2}{k}}n.
\end{align*}

\section{Bounds on $f_{1}\left(  n\right)  $}

Before stating the main result of this section, we shall recall two auxiliary
results whose proofs can be found in \cite{Nik05}. Given a graph $G=G\left(
n,m\right)  ,$ let%
\[
s\left(  G\right)  =\sum_{u\in V\left(  G\right)  }\left\vert d\left(
u\right)  -\frac{2m}{n}\right\vert .
\]

\begin{proposition}
For every graph $G=G\left(  n,m\right)  ,$%
\begin{equation}
\frac{s^{2}\left(  G\right)  }{2n^{2}\sqrt{2m}}\leq\mu_{1}\left(  G\right)
-\frac{2m}{n}\leq\sqrt{s\left(  G\right)  }, \label{prpin1}%
\end{equation}
and%
\begin{equation}
\mu_{n}\left(  G\right)  +\mu_{n}\left(  \overline{G}\right)  \leq
-1-\frac{s^{2}\left(  G\right)  }{n^{3}}. \label{prpin2}%
\end{equation}

\end{proposition}

Decreasing the constant $\sqrt{2}$ in (\ref{Nosin}) happened to be a
surprisingly challenging task for the author. The little progress that has
been made is given in the following theorem.

\begin{theorem}
\label{thf1} There exists $c\geq8\times10^{-7}$ such that%
\[
\mu_{1}\left(  G\right)  +\mu_{1}\left(  \overline{G}\right)  \leq\left(
\sqrt{2}-c\right)  n.
\]
for every graph $G$ of order $n$.
\end{theorem}

\begin{proof}
Assume the opposite: let $\varepsilon=8\times10^{-7}$ and let there exist a
graph $G$ of order $n$ such that
\[
\mu_{1}\left(  G\right)  +\mu_{1}\left(  \overline{G}\right)  >\left(
\sqrt{2}-\varepsilon\right)  n.
\]
Writing $A\left(  G\right)  $ for the adjacency matrix of $G,$ we have
\begin{equation}
\sum_{i=1}^{n}\mu_{i}^{2}\left(  G\right)  =tr\left(  A^{2}\left(  G\right)
\right)  =2e\left(  G\right)  , \label{basin}%
\end{equation}
implying that%
\[
\mu_{1}^{2}\left(  G\right)  +\mu_{n}^{2}\left(  G\right)  +\mu_{1}^{2}\left(
\overline{G}\right)  +\mu_{n}^{2}\left(  \overline{G}\right)  \leq2e\left(
G\right)  +2e\left(  \overline{G}\right)  <n^{2}.
\]
From%
\[
\mu_{1}^{2}\left(  G\right)  +\mu_{1}^{2}\left(  \overline{G}\right)
\geq\frac{1}{2}\left(  \mu_{1}\left(  G\right)  +\mu_{1}\left(  \overline
{G}\right)  \right)  ^{2}>\left(  1-\frac{\varepsilon}{\sqrt{2}}\right)
^{2}n^{2}>\left(  1-2\varepsilon\right)  n^{2}%
\]
we find that
\begin{equation}
\left\vert \mu_{n}\left(  G\right)  \right\vert +\left\vert \mu_{n}\left(
\overline{G}\right)  \right\vert \leq\sqrt{2\left(  \mu_{n}^{2}\left(
G\right)  +\mu_{n}^{2}\left(  \overline{G}\right)  \right)  }<\sqrt
{4\varepsilon}n^{2}, \label{in1}%
\end{equation}
and so, $\mu_{n}\left(  G\right)  +\mu_{n}\left(  \overline{G}\right)
>-\sqrt{4\varepsilon}n.$ We thus have $\sqrt{4\varepsilon}n^{4}\geq
s^{2}\left(  G\right)  .$ On the other hand, by (\ref{prpin1}) and in view of
$s\left(  G\right)  =s\left(  \overline{G}\right)  ,$ we see that
\[
\mu_{1}\left(  G\right)  +\mu_{1}\left(  \overline{G}\right)  \leq
n-1+2\sqrt{s\left(  G\right)  }<n+2\sqrt{s\left(  G\right)  },
\]
and, by (\ref{in1}), it follows that
\[
\left(  \sqrt{2}-\varepsilon\right)  n<n+2\left(  4\varepsilon\right)
^{1/8}n.
\]
Dividing by $n$, we obtain $\left(  \sqrt{2}-1\right)  <\varepsilon
+2^{5/4}\varepsilon^{1/8}$, a contradiction for $\varepsilon=8\times10^{-7}.$
\end{proof}

It is certain that the upper bound given by Theorem \ref{thf1} is far from the
best one. We shall give below a lower bound on $f_{1}\left(  n\right)  $ which
seems to tight.

For every $1\leq r<n$ the graph $G=K_{r}+\overline{K_{n-r}}$ (see, e.g.
\cite{HSF01}) satisfies
\[
\mu_{1}\left(  G\right)  +\mu_{1}\left(  \overline{G}\right)  =\frac{r-1}%
{2}+\sqrt{nr-\frac{3r^{2}+2r-1}{4}}+n-r-1=n-\frac{r+3}{2}+\sqrt{nr-\frac
{3r^{2}+2r-1}{4}}.
\]
The right-hand side of this inequality is increasing in $r$ for $0\leq
r\leq\left(  n-1\right)  /3$ and we find that
\[
f_{1}\left(  n\right)  >\frac{4n}{3}-2.
\]
This gives some evidence for the following conjecture.

\begin{conjecture}%
\[
f_{1}\left(  n\right)  =\frac{4n}{3}+O\left(  1\right)  .
\]

\end{conjecture}

We conclude this section with an improvement of the lower bound in
(\ref{Nosin}). Using the first of inequalities (\ref{prpin1}) we obtain
\begin{align*}
\mu_{1}\left(  G\right)  +\mu_{1}\left(  \overline{G}\right)   &  \geq
n-1+\frac{s^{2}\left(  G\right)  }{2n^{2}}\left(  \frac{1}{\sqrt{2e\left(
G\right)  }}+\frac{1}{\sqrt{2e\left(  \overline{G}\right)  }}\right)  \geq\\
&  \geq n-1+\sqrt{2}\frac{s^{2}\left(  G\right)  }{n^{3}}.
\end{align*}

\section{A class of graphs}

In this section we shall describe a class of graphs that give the right order
of $f_{2}\left(  G\right)  $ and, we believe, also of $f_{n}\left(  G\right)
.$

Let $n\geq4$. Partition $\left[  n\right]  $ in $4$ classes $A,B,C,D$ so that
$\left\vert A\right\vert \geq\left\vert B\right\vert \geq\left\vert
C\right\vert \geq\left\vert D\right\vert \geq\left\vert A\right\vert -1.$ Join
every two vertices inside $A$ and $D,$ join each vertex in $B$ to each vertex
in $A\cup C,$ join each vertex in $D$ to each vertex in $C.$ Write $G\left(
n\right)  $ for the resulting graph.

Note that if $n$ is divisible by $4,$ the sets $A,B,C,D$ have equal
cardinality and we see that $G\left(  n\right)  $ is isomorphic to its complement.

Our main goal to the end of this section is to estimate the eigenvalues of
$G\left(  n\right)  .$ Write $ch\left(  A\right)  $ for the characteristic
polynomial of a matrix $A.$ The following general theorem holds.

\begin{theorem}
\label{thch}Suppose $G$ is a graph and $V\left(  G\right)  =\cup_{i=1}%
^{k}V_{i}$ is a partition in sets of size $n$ such that

(i) for all $1\leq i\leq k,$ either $e\left(  V_{i}\right)  =\binom{n}{2}$ or
$e\left(  V_{i}\right)  =0$;

(ii) for all $1\leq i<j\leq k,$ either $e\left(  V_{i},V_{j}\right)  =n^{2}$
or $e\left(  V_{i},V_{j}\right)  =0.$

Let the sets $V_{1},...,V_{p}$ be independent and $V_{p+1},...,V_{k}$ induce a
complete graph. Then for the characteristic polynomial of the adjacency matrix
of $G$ we have%
\[
ch\left(  A\left(  G\right)  \right)  =x^{pn-p}\left(  1-x\right)  ^{\left(
k-p\right)  n-\left(  k-p\right)  }ch\left(  R\right)  ,
\]
where $R=\left(  r_{ij}\right)  $ is a $k\times k$ matrix such that
\[
r_{ij}=\left\{
\begin{array}
[c]{lll}%
0 & \text{if }i\neq j & \text{and }e\left(  V_{i},V_{j}\right)  =0\\
n & \text{if }i\neq j & \text{and }e\left(  V_{i},V_{j}\right)  =n^{2}\text{
}\\
0 & \text{if }i=j & \text{and }e\left(  V_{i}\right)  =0\\
n-1 & \text{if }i=j & \text{and }e\left(  V_{i}\right)  =\binom{n}{2}%
\end{array}
\right.  .
\]

\end{theorem}

The proof of this theorem is a straight exercise in determinants, so we shall
omit it.

If $n$ is divisible by $4,$ say $n=4k,$ by Theorem \ref{thch}, for the
characteristic polynomial of $A\left(  G\left(  n\right)  \right)  $ we have%
\[
ch\left(  A\left(  G\left(  n\right)  \right)  \right)  =x^{2k-2}\left(
1-x\right)  ^{2k-2}\left[
\begin{array}
[c]{cccc}%
k-1-x & k & 0 & 0\\
k & -x & k & 0\\
0 & k & -x & k\\
0 & 0 & k & k-1-x
\end{array}
\right]  .
\]
By straightforward calculations, setting $a=1-1/k$ and $y=x/k$, we see that
\begin{align*}
ch\left(  A\left(  G\left(  n\right)  \right)  \right)   &  =x^{2k-2}\left(
1-x\right)  ^{2k-2}\left[  \left(  a-y\right)  \left(  y^{2}\left(
a-y\right)  +2y-a\right)  -\left(  y^{2}-ay-1\right)  \right] \\
&  =x^{2k-2}\left(  1-x\right)  ^{2k-2}\left(  y^{2}-\left(  1+a\right)
y-\left(  1-a\right)  \right)  \left(  y^{2}+\left(  1-a\right)  y-\left(
a+1\right)  \right)  .
\end{align*}
Hence, we find that%
\begin{align*}
\mu_{2}\left(  G\right)   &  =-\frac{1}{2}+\sqrt{\frac{1}{4}+2\left\lfloor
\frac{n}{4}\right\rfloor ^{2}-\left\lfloor \frac{n}{4}\right\rfloor }\\
\mu_{n}\left(  G\right)   &  =-\frac{1}{2}-\sqrt{\frac{1}{4}+2\left\lfloor
\frac{n}{4}\right\rfloor ^{2}-\left\lfloor \frac{n}{4}\right\rfloor .}%
\end{align*}

If $n$ is not divisible by $4,$ we will give some tight estimates of $\mu
_{2}\left(  G\right)  $ and $\mu_{n}\left(  G\right)  .$ Notice first that
$G\left(  4\left\lfloor n/4\right\rfloor \right)  $ is an induced graph of
$G\left(  n\right)  $ which in turn is an induced graph of $G\left(
4\left\lceil n/4\right\rceil \right)  .$ Thus the adjacency matrix of
$G\left(  4\left\lfloor n/4\right\rfloor \right)  $ is a principal submatrix
of the adjacency matrix of $G\left(  n\right)  $ which in turn is a principal
submatrix of the adjacency matrix of $G\left(  4\left\lceil n/4\right\rceil
\right)  .$ Since the eigenvalues of a matrix and its principal matrices are
interlaced (\cite{HoJo88}, Theorem 4.3.15), we obtain%
\begin{align}
-\frac{1}{2}+\sqrt{\frac{1}{4}+2\left\lfloor \frac{n}{4}\right\rfloor
^{2}-\left\lfloor \frac{n}{4}\right\rfloor }  &  \leq\mu_{2}\left(  G\right)
\leq-\frac{1}{2}+\sqrt{\frac{1}{4}+2\left\lceil \frac{n}{4}\right\rceil
^{2}-\left\lceil \frac{n}{4}\right\rceil },\label{in3}\\
-\frac{1}{2}-\sqrt{\frac{1}{4}+2\left\lfloor \frac{n}{4}\right\rfloor
^{2}-\left\lfloor \frac{n}{4}\right\rfloor }  &  \leq\mu_{n}\left(  G\right)
\leq-\frac{1}{2}-\sqrt{\frac{1}{4}+2\left\lceil \frac{n}{4}\right\rceil
^{2}-\left\lceil \frac{n}{4}\right\rceil }. \label{in4}%
\end{align}

\section{The asymptotics of $f_{2}\left(  n\right)  $}

In this section we shall prove inequalities (\ref{mainin2}). From (\ref{in3})
we readily have%
\[
f_{2}\left(  n\right)  \geq-\frac{1}{2}+\sqrt{\frac{1}{4}+2\left\lfloor
\frac{n}{4}\right\rfloor ^{2}-\left\lfloor \frac{n}{4}\right\rfloor }%
>\frac{\sqrt{2}}{2}n-3,
\]
so all we need to prove is that $f_{2}\left(  n\right)  \leq n/\sqrt{2}$.

By (\ref{basin}) we have
\begin{equation}
\mu_{1}^{2}\left(  G\right)  +\mu_{2}^{2}\left(  G\right)  +\mu_{n}^{2}\left(
G\right)  +\mu_{1}^{2}\left(  \overline{G}\right)  +\mu_{2}^{2}\left(
\overline{G}\right)  +\mu_{n}^{2}\left(  \overline{G}\right)  \leq n\left(
n-1\right)  . \label{in2}%
\end{equation}
By Weyl's inequalities (\cite{HoJo88}, p. 181), for every graph $G$ of order
$n,$ we have
\[
\mu_{2}\left(  G\right)  +\mu_{n}\left(  \overline{G}\right)  \leq\mu
_{2}\left(  K_{n}\right)  =-1.
\]
Hence, using $\mu_{2}\geq0$ and $\mu_{n}\leq-1$ we obtain
\[
\mu_{2}^{2}\left(  G\right)  \leq\mu_{n}^{2}\left(  \overline{G}\right)
+2\mu_{n}\left(  \overline{G}\right)  +1<\mu_{n}^{2}\left(  \overline
{G}\right)  .
\]
Hence, from (\ref{in2}) and $\mu_{1}\left(  G\right)  +\mu_{1}\left(
\overline{G}\right)  \geq n-1,$ we find that
\[
\frac{\left(  n-1\right)  ^{2}}{2}+2\mu_{2}^{2}\left(  G\right)  +2\mu_{2}%
^{2}\left(  \overline{G}\right)  \leq\mu_{1}^{2}\left(  G\right)  +\mu_{2}%
^{2}\left(  G\right)  +\mu_{n}^{2}\left(  G\right)  +\mu_{1}^{2}\left(
\overline{G}\right)  +\mu_{2}^{2}\left(  \overline{G}\right)  +\mu_{n}%
^{2}\left(  \overline{G}\right)  \leq n\left(  n-1\right)  .
\]
After some algebra, we deduce that
\[
\mu_{2}\left(  G\right)  +\mu_{2}\left(  \overline{G}\right)  \leq\frac
{\sqrt{2}}{2}n,
\]
completing the proof of inequalities (\ref{mainin2}).

\section{Bounds on $f_{n}\left(  n\right)  $}

In this section we shall prove inequalities (\ref{mainin3}). From (\ref{in4}),
as above, we have%
\[
f_{n}\left(  n\right)  >\frac{\sqrt{2}}{2}n-3.
\]
We believe that, in fact, the following conjecture is true.

\begin{conjecture}%
\[
f_{n}\left(  G\right)  =\frac{\sqrt{2}n}{2}+O\left(  1\right)  .
\]

\end{conjecture}

However we can only prove that $f_{n}\left(  G\right)  <\left(  \sqrt
{3}/2\right)  n$ which is implied by the following theorem.

\begin{theorem}
For every graph $G$ of order $n,$
\[
\mu_{n}^{2}\left(  G\right)  +\mu_{n}^{2}\left(  \overline{G}\right)
\leq\frac{3}{8}n^{2}.
\]

\end{theorem}

\begin{proof}
Indeed, suppose $\left(  u_{1},...,u_{n}\right)  $ and $\left(  w_{1}%
,...,w_{n}\right)  $ are eigenvectors to $\mu_{n}\left(  G\right)  $ and
$\mu_{n}\left(  \overline{G}\right)  .$ Let
\[
U=\left\{  i:u_{i}>0\right\}  ,\text{ \ \ \ }W=\left\{  i:w_{i}>0\right\}  .
\]
Setting $V=\left[  n\right]  ,$ we clearly have $\mu_{n}^{2}\left(  G\right)
\leq E_{G}\left(  U,V\backslash U\right)  $ and $\mu_{n}^{2}\left(
\overline{G}\right)  \leq E_{\overline{G}}\left(  W,V\backslash W\right)  $.
Since $E_{G}\left(  U,V\backslash U\right)  \cap E_{\overline{G}}\left(
W,V\backslash W\right)  =\varnothing,$ we see that the graph
\[
G^{\prime}=\left(  V,E_{G}\left(  U,V\backslash U\right)  \cup E_{\overline
{G}}\left(  W,V\backslash W\right)  \right)
\]
is at most $4$-colorable and hence $G^{\prime}$ contains no $4$-cliques. By
Tur\'{a}n's theorem (e.g., see \cite{Bol98}), we obtain $e\left(  G^{\prime
}\right)  \leq\left(  3/8\right)  n^{2},$ completing the proof.
\end{proof}

\section{Bounds on $f_{k}\left(  n\right)  ,$ $2<k<n$}

In this section we shall give simple bounds on $f_{k}\left(  n\right)  $ for
$2<k<n$. Write for the Tur\'{a}n graph of order $n$ with $k$ classes. Recall
that $T_{k}\left(  n\right)  $ is a complete $k$-partite graph whose vertex
classes differ by at most $1$ in size. We assume that $k$ is fixed and $n$ is
large enough. Since $\mu_{k}\left(  T_{k}\left(  n\right)  \right)  =0$ and
$\mu_{n-k}\left(  T_{k}\left(  n\right)  \right)  \leq-\left\lfloor
n/k\right\rfloor $ for $n$ large, we immediately have
\begin{align*}
f_{k}\left(  n\right)   &  \geq\left\lfloor n/k\right\rfloor -1,\\
f_{n-k}\left(  n\right)   &  \geq\left\lfloor n/k\right\rfloor +1.
\end{align*}

We next turn to upper bounds on $f_{k}\left(  n\right)  .$

\begin{theorem}
For any fixed $k$ and any graph $G$ of sufficiently large order $n,$
\begin{equation}
\left\vert \mu_{k}\left(  G\right)  \right\vert +\left\vert \mu_{k}\left(
\overline{G}\right)  \right\vert <\sqrt{\frac{2}{k}}n. \label{in5}%
\end{equation}
and%
\begin{equation}
\left\vert \mu_{n-k}\left(  G\right)  \right\vert +\left\vert \mu_{n-k}\left(
\overline{G}\right)  \right\vert <\sqrt{\frac{2}{k}}n. \label{in6}%
\end{equation}

\end{theorem}

\begin{proof}
Set $e\left(  G\right)  =m.$ Our first goal is to prove that $\left\vert
\mu_{k}\left(  G\right)  \right\vert \leq\sqrt{2e\left(  G\right)  /k}.$ If
$\mu_{k}\left(  G\right)  \geq0,$ we have in view of (\ref{basin})
\[
k\mu_{k}^{2}\left(  G\right)  \leq\sum_{i=1}^{n}\mu_{i}^{2}\left(  G\right)
=2m.
\]
If $\mu_{k}\left(  G\right)  <0$ and $\left\vert \mu_{k}\left(  G\right)
\right\vert >\sqrt{2m/k}$ then
\[
\sum_{i=1}^{n}\mu_{i}^{2}\left(  G\right)  \geq\left(  n-k\right)  \mu_{k}%
^{2}\left(  G\right)  >2m\frac{n-k}{k}>2m,
\]
a contradiction. Hence, $\left\vert \mu_{k}\left(  G\right)  \right\vert
\leq\sqrt{2e\left(  G\right)  /k},$ and, by symmetry, $\left\vert \mu
_{k}\left(  \overline{G}\right)  \right\vert \leq\sqrt{2e\left(  \overline
{G}\right)  /k}.$ Now
\[
\left\vert \mu_{k}\left(  G\right)  \right\vert +\left\vert \mu_{k}\left(
\overline{G}\right)  \right\vert \leq\sqrt{2e\left(  G\right)  /k}%
+\sqrt{2e\left(  \overline{G}\right)  /k}\leq\sqrt{\frac{2}{k}n\left(
n-1\right)  }<\sqrt{\frac{2}{k}}n,
\]
proving inequality (\ref{in5}). The proof of inequality (\ref{in6}) goes along
the same lines, so we will omit it.
\end{proof}

\end{document}